\documentclass[letterpaper,12pt]{article}
\usepackage{tabularx} 
\usepackage{amsmath}  
\usepackage{graphicx} 
\usepackage[margin=1in,letterpaper]{geometry} 
\usepackage{cite} 
\usepackage[final]{hyperref} 
\hypersetup{
	colorlinks=true,       
	linkcolor=blue,        
	citecolor=blue,        
	filecolor=magenta,     
	urlcolor=blue         
}
\usepackage{blindtext}

\usepackage{amssymb}
\usepackage{titling}
\newcommand{\subtitle}[1]{%
  \posttitle{%
    \par\end{center}
    \begin{center}\Large#1\end{center}
    \vskip0.5em}%
}

\newtheorem{proposition}{Proposition}

\begin{document}

\title{Dynamics of SIR model with heterogeneous response to intervention policy}
\author{Dmitrii Rachinskii, Samiha Rouf}
\date{}
\maketitle


\begin{abstract}
	We study dynamics of a variant of the SIR the model, where we assume that individuals respond differently to dynamics of the epidemic. Their heterogeneous response is modeled by the Preisach hysteresis operator. 
	The degree of heterogeneity of the response is measured by the variance $\sigma^2$ of the corresponding distribution (the Preisach density function).
	The proposed model has a continuum of endemic equilibrium states characterized by different proportions of susceptible, infected and recovered populations. We consider how the limit point of the epidemic trajectory 
	and the infection peak along this trajectory depend on $\sigma$. The heterogeneous model is compared to the the switched model with an ideally uniform response corresponding to $\sigma=0$.
\end{abstract}

\bigskip
{\it 2010 Mathematics Subject Classification}: {\ 92D30, 92D25, 47J40 }

\bigskip
{\it Keywords}: {\ Switched transmission rate, heterogeneous response, switched system, Preisach operator, hysteresis, continuum of equilibrium states, stability}


\section{Introduction}

At the threat of epidemics, health and government authorities can intervene by 
raising awareness in the population about the current severity of the epidemic, providing access to 
tests, vaccination and medicines, working with school authorities, using media and administrative pressure, etc (Center of Disease Control and Prevention, 2021).
During the Covid-19 epidemic, 
disease prevention measures have been imposed 
on an unprecedented scale 
and included 
massive quarantine and social distancing measures, 
 business restrictions, gathering and travel limitations, transition to online teaching at schools and colleges
 and stay-at-home and  shelter-in-place orders (Gostin and Wiley, 2020; Marquioni and de Aguiar, 2020). 
One major objective of such measures 
is flattening the curve, i.e.\ slowing down the spread of the epidemic in order to keep the number of active disease cases at a manageable level 
dictated by the  capacity of the health care system (Matrajt and Leung, 2020). 
However, most community isolation and business restriction measures can be introduced only for a limited period of time by economic and social reasons (Fairlie, 2020). Due to these constraints, many of the intervention protocols, which have been implemented by the health authorities in order to contain the Covid-19 epidemic (Wilder-Smith et al., 2020), can be thought of, at least simplistically, as 
threshold based.
As such, an intervention starts when a certain variable such as the number of daily new infections, the percentage of the occupied hospital beds or the basic reproduction number $R_0$ reaches a critical threshold value set by the health or government authority (Department of Health and Human Services, Nebraska, 2020; DeBenedetto and Ruiz, 2021); the intervention is revoked when this variable drops below the level deemed safer.
If the thresholds at which the intervention begins and ends are different, then hysteresis effect is present.

The effectiveness of the prevention measures can depend on multiple factors and a complex interplay between them.
Indeed, significantly different Covid-19 disease statistics have been reported by countries, which implemented seemingly similar intervention strategies. One such factor is the response of the population to the intervention policies and, in particular, the degree of uniformity of the response. 
The ability and willingness of an individual to receive immunization 
or follow community isolation policies
depends on the perceived risk of contracting the disease, risk of possible complications, level of trust to the authorities in the community, personal beliefs, etc (Guidry et al., 2021; Lazarus et al., 2021). 
Social interaction, reinforcement and imitation  
can lead to a hysteresis effect in adopting a `healthy behavior'
(Su et al., 2017). For example,
	in the presence of imperfect vaccine, 
hysteresis loops of vaccination rate arise with respect to changes in the perceived cost of vaccination (Chen and Fu, 2019). 
This cost changes dynamically 
as
individuals revisit their vaccination decision in response to dynamics of the epidemic and through a social learning processes under peer influence.
One of the findings in (Chen and Fu, 2019) is that hysteresis becomes more pronounced 
with increasing heterogeneity of the population.

The risk of contracting the disease and the risks of possible complications vary with age, health condition, living circumstances and profession. On the other hand, interventions of the health authorities and administrative measures at the level of a county or state can vary in scale depending on the availability of resources, the local economic situation and other factors (Emanuel et al., 2020). All these variations lead to the heterogeneity of the response of the population to the advent of an epidemic.

In this paper, we attempt to model the effect of the heterogeneity of the response on the epidemic trajectory
in the presence of hysteresis.
To this end, we use a variant of the SIR model where the transmission rate depends on dynamics of the infected population.
As a starting point, we adopt the approach of Chladn\'a et al. (2020) to modeling the uniformly homogeneous switched response of a population to the varying number of infected individuals by a two-threshold two-state relay operator (see Section 2.2).
In this model, it is assumed that the health authority implements a two-threshold intervention policy
whereby the intervention starts when the number of infected individuals exceeds a critical threshold value and stops whenever this number drops below a different (lower) threshold. 
The objective of the two-threshold policy is to navigate the system to the endemic equilibrium simultaneously keeping the number of infected individuals in check and not committing to continuous intervention. If the response is ideally uniform, then  prevention measures are assumed to
 translate immediately to reduced values of the transmission rate and $R_0$ during the intervention.
Next, to reflect the heterogeneity of the response, 
the population is divided into multiple subpopulations,
each characterized by a different pair of switching thresholds. In order to keep the model relatively simple, we apply averaging under further simplifying assumptions. The main simplification is that perfect mixing of the population is assumed. This leads to a differential model with just two variables, $S$ and $I$, but with a complex operator relationship between the transmission rate and the density of the infected population $I$.  As such, this operator relationship, known as the Preisach hysteresis operator (see, for example, Krasnosel'skii et al., 1989; Mayergoyz 1993; Visintin 1994; Brokate et al., 1996; Krej\v{c}\'\i, 1996) accounts for the heterogeneity of the response. This approach is similar to the one in Pimenov et al., 2012.
Heterogeneity of intervention policies can be modeled in a similar fashion
(see Section 3.2).

The heterogeneous model proposed below has a continuum of endemic equilibrium states characterized by different proportions of the susceptible, infected and recovered populations. 
We consider the convergence of the epidemic trajectory to this continuum.
Further, using the variance $\sigma^2$ of the distribution of the susceptible population over the set of switching thresholds
as a measure of the degree of heterogeneity of the response, we 
consider how parameters of the epidemic trajectory depend on $\sigma$. These parameters include
the peak of infection along the epidemic trajectory
and the proportions of the susceptible, infected and recovered populations at the endemic state where the trajectory converges to.

The paper is organized as follows. Models with an ideally homogeneous and heterogeneous response are discussed in Sections 3 and 4, respectively. Numerical results are discussed in Section 4. Proofs are presented in the Appendix.



\section{System with switched transmission coefficient}

{\bf 2.1.} The systems considered below adapt the standard scaled SIR model
\begin{equation}\label{basemodel}
\begin{aligned}
    \dot I &= \beta SI - (\gamma + \mu)I,\\
    \dot S &= -\beta SI - \mu S + \mu,
\end{aligned}
\end{equation}
where $I$ and $S$ are the densities of the infected and susceptible populations, respectively; $\beta$ is the transmission coefficient; $\gamma$ is the recovery rate; and, $\mu$ is the departure rate
due to the disease unrelated death, emigration etc. We assume a constant total population scaled to unity, hence the density of the recovered population $R=1-I-S$ can be removed from the system. 
The domain
$I, S\ge 0$, $ I+S\le 1$ is positively invariant for system \eqref{basemodel}.

Since the transmission coefficient is not measured in practice, we can re-scale time in (\ref{basemodel}) to obtain the normalized system  
\begin{equation}\label{r0model}
\begin{aligned}
    \dot I &= R_0 SI - I, \\
    \dot S &= - R_0 SI - \rho S + \rho,
\end{aligned}
\end{equation}
with the dimensionless parameters 
$$
R_0 = \frac{\beta}{\gamma+\mu}, \qquad \rho = \frac{\mu}{\gamma+\mu}<1.
$$
If the {\em basic reproduction number} satisfies $R_0<1$, then the {\em infection free} equilibrium $(I_*,S_*)=(0,1)$
is globally stable in the closed positive quadrant. On the other hand, if $R_0<1$, then
the infection free equilibrium is a saddle, and the  {\em endemic} equilibrium $(I^*,S^*)$ defined by
\begin{equation}\label{equil}
I^* = \left(1-\frac{1}{R_0}\right)\rho, 
\qquad S^* = \frac{1}{R_0}
\end{equation}
is globally stable in the open positive quadrant.
If 
\begin{equation}\label{focus}
\rho(R_0)^2 < 4(R_0 - 1), 
\end{equation}
then the endemic equilibrium is of focus type.

\begin{figure}[h]
\centering
\includegraphics[width = 0.5\linewidth]{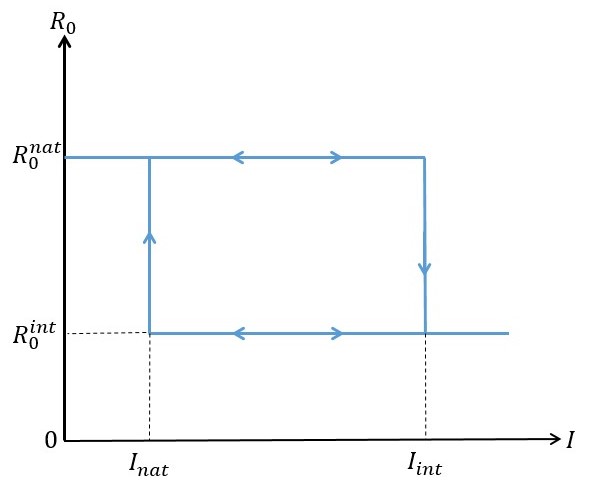}
\caption{Switching rule \eqref{relay1eq}. The vertical segments correspond to instantaneous transitions
of the basic reproduction number from the value $R_0^{nat}$ to $R_0^{int}$ and back. \label{relay1fig}}
\end{figure}

\medskip
{\bf 2.2.} We consider a switched system with two flows 
of the form (\ref{basemodel}) with two different values of the basic reproduction number,
$R_0=R_0^{nat}$ and $R_0=R_0^{int}$, where 
\begin{equation}\label{RR}
R_0^{nat} > R_0^{int}>1; 
\end{equation}
the parameter $\rho$ is the same for both flows. 
A switch occurs when the density of the infected population reaches certain thresholds, $I_{int}$ and $I_{nat}$,
with 
\begin{equation}\label{II}
0<I_{nat} < I_{int}<1. 
\end{equation}
We postulate that the basic reproduction number instantaneously switches from 
the value $R_0^{nat}$ to the value $R_0^{int}$ when the variable $I=I(t)$ reaches the upper threshold value $I_{int}$.
On the other hand, $R_0$ switches back from the value $R_0^{int}$ to the value $R_0^{nat}$ as $I$ reaches the lower
threshold value $I_{nat}$, see Figure \ref{relay1fig}.
%
%
%
This two-threshold switching rule can be formalized as follows:
\begin{equation}\label{relay1eq}
R_0(t)=
\left\{
\begin{array}{ll}
R_0^{nat} & \text{if either $
	I(\tau)<I_{int}$ for all $\tau\in[0,t]$}\\
& \text{or there is $t_1\in [0,t]$ such that $I(t_1)\le I_{nat}$}\\
& \text{and $I(\tau)<I_{int}$ for all $\tau\in (t_1, t]$;}\\
R_0^{int} & \text{if there is $t_1\in [0,t]$ such that $I(t_1)\ge I_{int}$}\\
& \text{and $I(\tau)>I_{nat}$ for all $\tau\in (t_1, t]$,}
\end{array}\right.
	\end{equation}
where for simplicity we assume that initially $I(0)<I_{int}$, which is sufficient for our purposes.
Hence, the switched system is defined by equations \eqref{r0model} in which $R_0=R_0(t)$
is given by formula \eqref{relay1eq}. This system will be considered in the positively invariant domain 
$I, S\ge 0$, $ I+S\le 1$.

According to the interpretation discussed in the Introduction, switched system \eqref{r0model}, \eqref{relay1eq} 
models the dynamics of the epidemic under the assumption that the health authorities implement the two-threshold intervention policy --- 
an intervention begins when the number of active cases exceeds the threshold value $I_{nat}$ 
and is revoked when the number of active cases drops 
below the lower 
threshold value $I_{nat}$.
It is assumed that the intervention quickly translates into the reduction of the basic reproduction number; 
once the intervention stops,
$R_0$ returns to the larger value. 

\medskip 
{\bf 2.3.} 
Switched system \eqref{r0model}, \eqref{relay1eq} has the infection free equilibrium $(I_*,S_*)=(0,1)$
of saddle type.
In addition, depending on the relative positions of the the threshold lines $I=I_{nat}$, $I=I_{int}$
and points 
\begin{equation}\label{EE}
E_{nat}=(I^*_{nat},S^{*}_{nat})=\left(\left(1-\frac{1}{R_0^{nat}}\right)\rho,\frac1{R_0^{nat}}\right),\quad E_{int} =
(I^*_{int},S^{*}_{int})=\left(\left(1-\frac{1}{R_0^{int}}\right)\rho,\frac1{R_0^{int}}\right),
\end{equation}
the switched system can additional have
either one or two stable endemic equilibrium states;
here $E_{nat}$ is the endemic equilibrium of system \eqref{r0model} with constant $R_0=R_0^{nat}$,
while $E_{int}$ is the endemic equilibrium of \eqref{r0model} with constant $R_0=R_0^{int}$ (cf.~\eqref{equil}).
The following statement follows directly from the switching rule.

\begin{proposition}\label{p0}
Switched system \eqref{r0model}, \eqref{relay1eq} has 
\begin{itemize}
\item[--] no endemic equilibrium states if $I^*_{int}<I_{nat}<I_{int}<I^*_{nat}$;
\item[--] one endemic equilibrium state $E_{nat}$ if 
$I^*_{int}<I_{nat}$ and $I^*_{nat}<I_{int}$;
\item[--] one endemic equilibrium state $E_{int}$ if 
$I_{nat}<I^*_{int}$ and $I_{int}<I^*_{nat}$;
\item[--] two endemic equilibrium states $E_{nat}$ and $E_{int}$ if $I_{nat}<I^*_{int}<I^*_{nat}<I_{int}$.
\end{itemize}
Each endemic equilibrium is stable. 
\end{proposition}

In the last case of the above alternative, which is characterized by the bi-stability,
the equilibrium $E_{int}$ has a higher population of susceptible individuals and lower populations of infected and recovered individuals than the equilibrium $E_{nat}$.

\begin{figure}[ht!] \label{relayfig}
    \centering
    \text{(a)}\includegraphics[width = 0.32\linewidth]{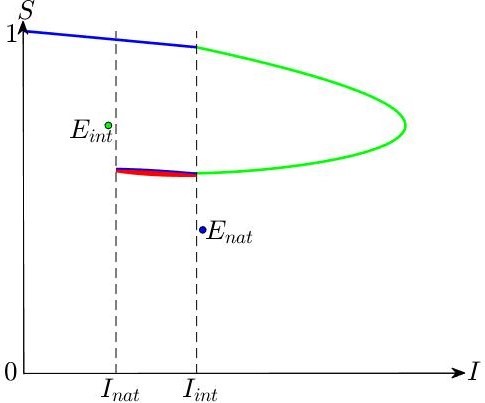}\hspace{1in}
    \text{(b)}\includegraphics[width = 0.32\linewidth]{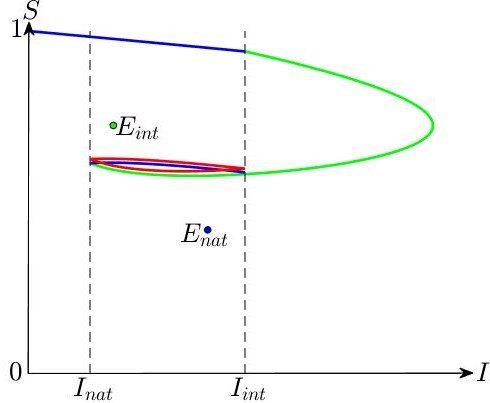}\hspace{1in}
    \text{(c)}\includegraphics[width = 0.32\linewidth]{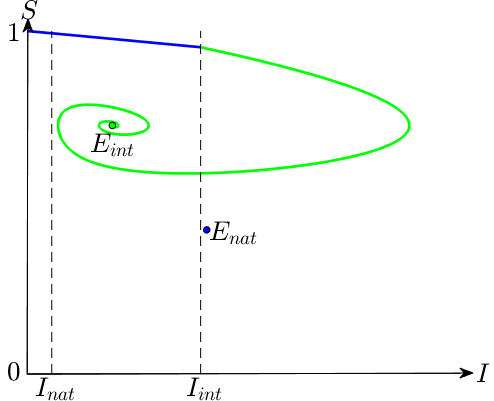}\hspace{1in}
   \text{(d)} \includegraphics[width = 0.32\linewidth]{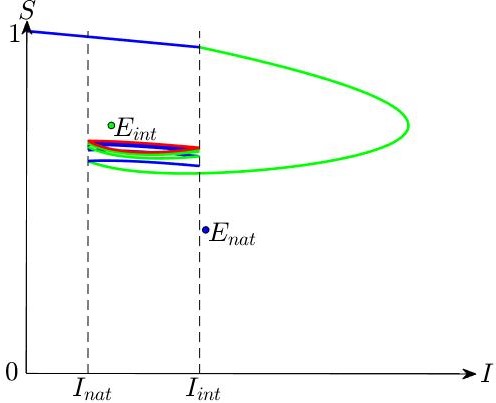}\hspace{1in}
   \text{(e)} \includegraphics[width = 0.32\linewidth]{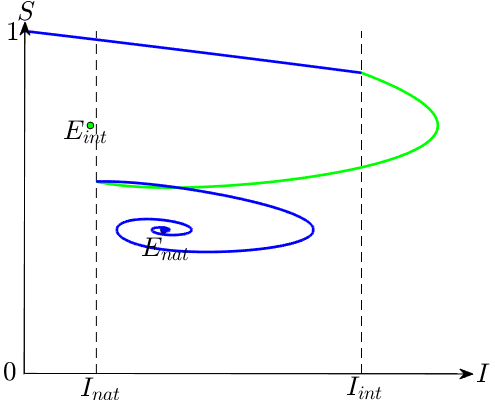}\hspace{1in}
   \text{(f)} \includegraphics[width = 0.32\linewidth]{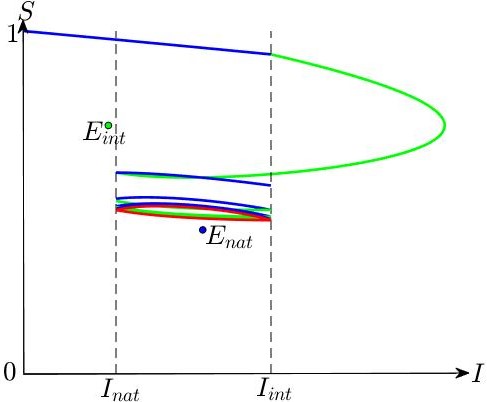}\hspace{1in}
    \text{(g)}\includegraphics[width = 0.32\linewidth]{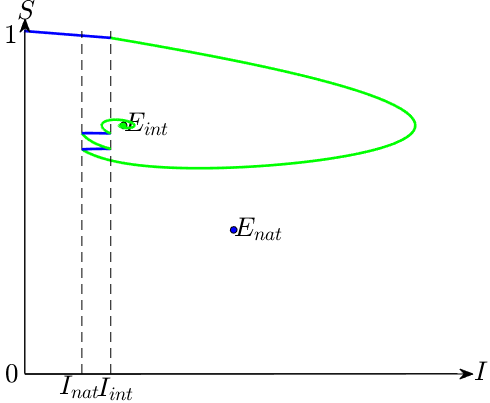}\hspace{1in}
    \text{(h)}\includegraphics[width = 0.32\linewidth]{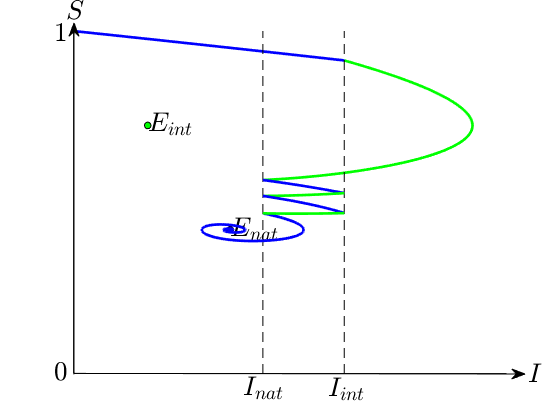}
    \caption{Convergence of a trajectory of switched system \eqref{r0model}, \eqref{relay1eq} to either an endemic equilibrium state or to a periodic orbit depending on the placement of the thresholds $I_{nat}$ and $I_{int}$. Blue segments of the trajectory correspond to $R_0(t)=R_0^{nat}$; green segments correspond to $R_0(t)=R_0^{int}$.
    Periodic attractor is shown by red. The parameters are $\rho = 0.05$, $R_0^{nat} = 2.38$, $R_0^{int} = 1.38$, $E_{int}= (0.013,0.73)$, $E_{nat}=(0.029,0.42)$, $I(0)=0.001$, $S(0)=0.999$.}
\end{figure}

\begin{figure}[ht!] 
    \centering
    \includegraphics[width = 0.7\linewidth]{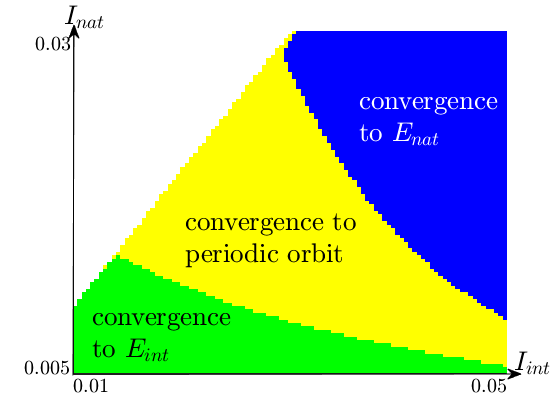}
    \caption{The color map indicates the convergence of a trajectory of switched system \eqref{r0model}, \eqref{relay1eq} to either $E_{int}$ or $E_{nat}$ or a periodic orbit depending on the threshold values $I_{int}$,$I_{nat}$ where $I_{int}>I_{nat}$. The parameters are $\rho = 0.05$, $R_0^{nat} = 2.38$, $R_0^{int} = 1.38$, $E_{int}= (0.013,0.73)$, $E_{nat}=(0.029,0.42)$, $I(0)=0.001$, $S(0)=0.999$.
    }\label{colormap}
\end{figure}

We are primarily interested in the trajectories that start near the infection free equilibrium $(I_*,S_*)=(0,1)$.
Figure \ref{relayfig} presents typical behaviors of such a trajectory depending on the placement of the thresholds $I=I_{nat}, I_{int}$ relative to the points $E_{nat}, E_{int}$. If $I^*_{int}\le I_{nat}<I_{int}\le I^*_{nat}$, then the trajectory converges to a periodic orbit as shown in Figures \ref{relayfig}(a,b,d,f). This scenario is robust. 
In particular, it persists  under small variations of the threshold values such that $I_{nat}<I_{int}^*$  and/or $I_{int}>I^*_{nat}$, see Figures \ref{relayfig}(d,f), respectively. However, for larger values of the upper threshold $I_{int}$ (with the other threshold fixed), we observe the convergence to the endemic equilibrium $E_{nat}$ as in Figure \ref{relayfig}(e). On the other hand, decreasing the lower threshold $I_{nat}$ results in the convergence to the endemic equilibrium $E_{int}$, see Figure \ref{relayfig}(c).


These observations are in agreement with the results proved for system \eqref{r0model}, \eqref{relay1eq} in the case  $R_0^{nat}>1>R_0^{int}$ in Chladn\'a et al. (2020), where it was shown that each trajectory converges either to a periodic orbit or to the endemic equilibrium.
Here we observe similar attractors in the case \eqref{RR}. In particular,
for a periodic orbit with $I(t)$ oscillating between the values $I_{nat}$ and $I_{int}$,
the point $(I(t),R_0(t))$ moves clockwise along the {rectangular hysteresis loop} shown in Figure~\ref{relay1fig}. 
 


\section{
Model with heterogeneous transmission}

{\bf 3.1.}
Let us recall the definition of the 
 so-called {\em non-ideal relay operator}, which is
 also known as {a rectangular hysteresis loop} or a lazy switch (see e.g.~Visintin 1994).
 This operator is characterized by two scalar parameters $\alpha_1$ and $\alpha_2$, the threshold values, with $\alpha_1<\alpha_2$; we will use the notation $\alpha=(\alpha_1,\alpha_2)$. The {\em input} of the relay is an arbitrary continuous function of time, $I: \mathbb{R}_+\to\mathbb{R}$. The {\em state} of the relay, denoted $r_{\alpha}(t)$, equals either $0$ or $1$ at any moment $t\in \mathbb{R}_+$.
 More specifically,
 given any continuous input $I: \mathbb{R}_+\to \mathbb{R}$ and an initial value of the state, $r_{\alpha}(0)=r_{\alpha}^0$,
which satisfies the constraints
 \begin{equation}\label{v1}
 r_{\alpha}^0\in\{0,1\} \quad \text{if} \quad \alpha_1<I(0)< \alpha_2;
 \end{equation}
 \begin{equation}\label{v2} 
 r_{\alpha}^0=0 \quad \text{if} \quad I(0)\le \alpha_1;\qquad r_{\alpha}^0=1 \quad \text{if} \quad I(0)\ge\alpha_2,
 \end{equation}
 the state of the relay at the future moments $t>0$ is defined by the relations
 \begin{equation}\label{relay'}
 r_{\alpha}(t)=
 \left\{
 \begin{array}{cl}
 0 & \text{if there is $t_1\in [0,t]$ such that $I(t_1)\le \alpha_1$}\\
 & \text{and $I(\tau)<\alpha_2$ for all $\tau\in (t_1, t]$;}\\
 1 & \text{if there is $t_1\in [0,t]$ such that $I(t_1)\ge \alpha_2$}\\
 & \text{and $I(\tau)>\alpha_1$ for all $\tau\in (t_1, t]$};\\
 r_{\alpha}(0) & \text{if $\alpha_1<I(\tau)<\alpha_2$ for all $\tau\in[0,t]$.}
 \end{array}\right.
 \end{equation}
 This function, 
 which depends both on the input $I(t)$ $(t\geq 0)$ and the initial state  $r_{\alpha}(0)=r_{\alpha}^0$  of the relay, will be denoted by
 \begin{equation}\label{re}
 r_{\alpha}(t)=({\mathcal R}_{\alpha}[r_{\alpha}^0]I)(t),\qquad t\geq 0.
 \end{equation}
 By definition \eqref{relay'} of this input-to-state map \eqref{re}, the state satisfies the constraints
 \begin{equation}\label{compatibility}
 r_{\alpha}(t)=1 \quad \text{whenever} \quad I(t)\geq
 \alpha_2; \quad r_{\alpha}(t)=0 \quad \text{whenever} \quad I(t)\leq \alpha_1
 \end{equation}
 at all times. Further,
 the function \eqref{relay'} has at most a finite number of jumps between the values $0$ and $1$
 on any finite time interval $t_{0}\leq t\leq t_{1}$. 
 
 

Using formula \eqref{relay'} and notation \eqref{re}, equation \eqref{relay1eq} for the basic reproduction number can be written equivalently as
\begin{equation}\label{rev}
R_0(t) = R_0^{nat} - \hat q \cdot (\mathcal{R}_{\alpha}[r_{\alpha}^0]I)(t),
\end{equation}
where $\alpha=(I_{nat},I_{int})$, 
$\hat q=R_0^{nat}-R_0^{int}$ and $r_\alpha^0=0$.

\medskip\noindent
{\bf 3.2.} Now, we consider a model, in which several responses of the form \eqref{rev}, with different thresholds 
$\alpha$, are combined 
because different groups of individuals respond differently to the advent and dynamics of an epidemic.

Individuals can respond differently to 
the interventions of the health authorities.
In particular, the ability and willingness to 
follow the recommendations of the health authority
can vary significantly from one to another group of individuals for the same level of threat of contracting the disease.
Multiple factors are at play such as the occupation, age, living environment and health condition of an individual, to mention a few.
In order to account for the heterogeneity of the individual response, let us divide the susceptible population into non-intersecting sub populations $S_\alpha$ parameterized by points $\alpha$ of a subset $\Pi\subset \{\alpha=(\alpha_1,\alpha_2): \alpha_1<\alpha_2\}$ of the $\alpha$-plane, assuming a homogeneous response within each group $S_\alpha$. As a simplification, let us assume at this point that the response of the infected population 
is heterogeneous and uniform, hence the transmission coefficient is determined by the behavior of the susceptible individuals only. 
Further, assume that 
the basic reproduction number for the sub population $S_\alpha$ is given by \eqref{rev} with 
$\hat q=\hat q(\alpha)$.
Then, the average basic reproduction number for the entire population at time $t$ equals
\begin{equation}\label{rev''}
R_0(t)= R_0^{nat} 
-\iint_{\Pi} \hat q(\alpha)\,
({\mathcal R}_{\alpha}[r_{\alpha}^0]I)(t)\, d F(\alpha), \qquad t\in\mathbb{R}_+,
\end{equation}
where the probability measure $F$ describes the distribution of the susceptible {population} $S$ over the index set $\Pi$ (the set of threshold pairs). Finally, we assume for simplicity that this distribution is independent of time,
i.e.~the probability measure $F$ does not change with variations of $I$. In this case, 
the mapping of the space of continuous
inputs $I:\mathbb{R}_+\to \mathbb{R}$ to the space of outputs $R_0:\mathbb{R}_+\to \mathbb{R}$
defined by \eqref{rev''} is known as the Preisach operator (see e.g.~Krasnosel'skii et al.~1983).

Under the above assumptions, the dynamics of the epidemic is modeled by system \eqref{r0model}
where $R_0(t)$ is related to the $I(t)$ by the Preisach operator \eqref{rev''}, which accounts 
for the heterogeneity of the transmission coefficient. 

A similar system results from the assumption 
that the health authorities 
have multiple intervention policies (numbered $n=1,\ldots,N$) 
in place, each decreasing the transmission coefficient by a certain amount $ \Delta \beta_n$ while the intervention is implemented.
The authorities aim to provide an adaptive response, which is adequate to the severity of the epidemic.
Let us suppose that each  intervention policy 
is guided by the two-threshold start/stop rule, such as in \eqref{rev}, associated with a particular pair of thresholds $\alpha^n=(\alpha_1^n,\alpha_2^n)$. 
The response of the population is assumed homogeneous.
Under these assumptions, the basic reproduction number of system \eqref{r0model} is given by
\begin{equation}\label{rev'}
R_0(t)= R_0^{nat} 
- \sum_{n=0}^{N-1} \hat q_n\cdot
({\mathcal R}_{\alpha^n}[r_{\alpha^n}^0]I)(t),
\end{equation}
where $\hat q_n=\Delta\beta_n/(\gamma+\mu)$;
hence, $R_0$ is set to change at multiple thresholds $\alpha_1^n, \alpha_2^n$.
Operator \eqref{rev'}, known as the discrete Preisach model, is a particular case of \eqref{rev''} with an atomic measure $F$. 

Below we consider absolutely continuous measures $F$. 
The corresponding operator \eqref{rev''}, which is called the continuous Preisach model,
can be written in the equivalent form
\begin{equation}\label{pre}
R_0(t)=R_0^{nat} - (R_0^{nat}-R_0^{int})\iint_{\Pi} q(\alpha)\, \bigl({\mathcal R}_{\alpha}[r_{\alpha}^0]I\bigr)(t) \,d\alpha_1 d\alpha_2,\qquad t\ge 0,
\end{equation}
where $0<R_0^{int}<R_0^{nat}$ and $q=q(\alpha): \Pi\to \mathbb{R}_+$ is a strictly positive probability measure density, i.e.
\begin{equation}\label{q}
\iint_\Pi q(\alpha)\,d\alpha_1d\alpha_2=1.
\end{equation}
This operator can be approximated by discrete operators \eqref{rev'}. We will assume that $q$ is bounded.

\medskip\noindent{{\bf 3.3.}}
Trajectories of system \eqref{r0model} coupled with the operator relationship \eqref{pre} lie in the infinite-dimensional phase space $\mathfrak U$ of triplets $(I,S,r^0)$, where
the measurable function $r^0=r_{\alpha}^0:\Pi\to \{0,1\}$ of the variable $\alpha=(\alpha_1,\alpha_2)$
describes
the states of the relays ${\mathcal R}_{\alpha}$ at the initial moment 
(see Appendix 6.1 for details).
Slightly abusing the notation, we will 
also refer to the two-dimensional curve  $(I(t),S(t))$ as a trajectory, omitting the component \eqref{re} in the state space of the Preisach operator.

Let us consider equilibrium states of system \eqref{r0model}, \eqref{pre}. The components $I$, $S$, $r^0$ of the solution and the basic reproduction number $R_0$ at an equilibrium state are constant, and $R_0$ is related to the function $r^0=r_{\alpha}^0:\Pi\to \{0,1\}$ by the equation
\begin{equation}\label{preR0}
R_0=R_0^{nat} - (R_0^{nat}-R_0^{int})\iint_{\Pi} q(\alpha) r_{\alpha}^0 \,d\alpha_1 d\alpha_2.
\end{equation}
In what follows,
\[
\Pi=\{\alpha=(\alpha_1,\alpha_2): 0\le \alpha_1<\alpha_2\le1\}.
\]
Due to this assumption and the compatibility constraint \eqref{compatibility}, the inclusion $\alpha\in\Pi$ implies that all the relays are in state $r_\alpha=0$ when $I=0$. Therefore,
system \eqref{r0model}, \eqref{pre} has a 
unique {\em infection free} equilibrium state $E_*=(I_*,S_*)=(0,1)$,
in which the function $r^0=r_\alpha^0: \Pi\to\{0,1\}$ is the identical zero and the basic reproduction number equals $R_0=R_0^{nat}$ according to \eqref{preR0}; it is the same equilibrium state as the switched system 
\eqref{r0model}, \eqref{relay1eq} has.

On the other hand,
system \eqref{r0model}, \eqref{pre} also has a {\em continuum} of {\em endemic} equilibrium states
\begin{equation}\label{endemic}
E_{\theta}=(I^*_{\theta},S^{*}_{\theta})=\left(\left(1-\frac{1}{R_0^{\theta}}\right)\rho,\frac1{R_0^{\theta}}\right), \qquad 0\le\theta\le 1,
\end{equation}
(cf.~\eqref{EE}), which are characterized by different proportions of the infected, susceptible and recovered populations
and different values $R_0^\theta$ of the basic reproduction number. More specifically, 
the following statement holds.

\begin{proposition}\label{p1}
Let relations \eqref{RR} and \eqref{q}
hold.
Then, endemic equilibrium states \eqref{endemic} of system \eqref{r0model}, \eqref{pre} form a non-degenerate line segment 
$E_0 E_1$ which belongs to the interior of the line segment $E_{int}E_{nat}$
connecting the 
points \eqref{EE}.
\end{proposition}

The proof is presented in Appendix 6.2.
\begin{figure}[ht!]
    \centering
    \includegraphics[width = 0.8\linewidth]{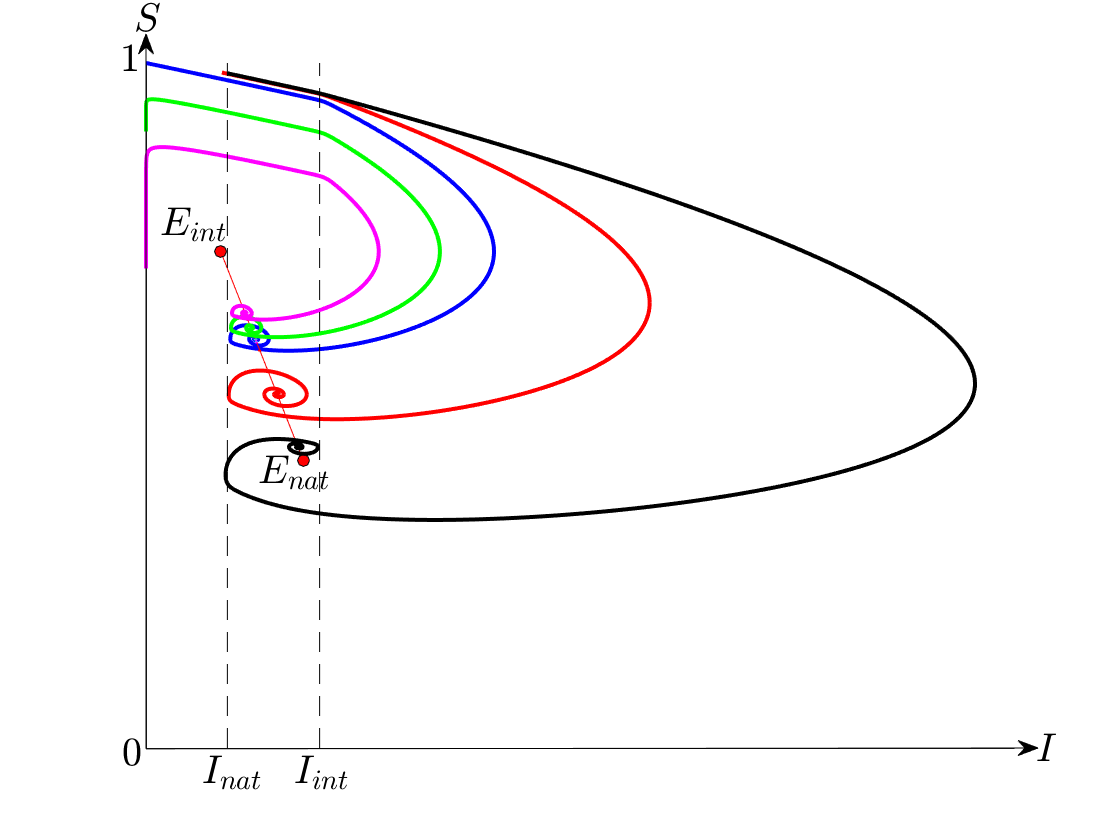}
    \caption{The limit point for different trajectories of system \eqref{r0model}, \eqref{pre}. Each trajectory converges to an equilibrium state, which belongs to the segment $E_{int}E_{nat}$. The density function of the Presiach operator is defined by \eqref{density}.
    The parameters are $\rho = 0.05$, $R_0^{nat} = 2.38$, $R_0^{int} = 1.38$, $E_{int}= (0.013,0.73)$, $E_{nat}=(0.029,0.42)$, $I_{nat}=0.015$, $I_{int}=0.032$, $\sigma=0.01$.}
    \label{prop2fig}
\end{figure}

\medskip\noindent{{\bf 3.4.} }
Let us consider a family of measure densities $q(\cdot;\sigma):\Pi\to\mathbb{R}_+$
and the corresponding probability measures on $\Pi$ depending on a parameter $\sigma>0$.
Assume that these measures converge to a Dirac measure concentrated at a point $\alpha^*\in \Pi$
 as $\sigma\to 0$, i.e.
\begin{equation}\label{sigma}
q(\alpha;\sigma) \to \delta(\alpha-\alpha^*)\qquad \text{as}\qquad \sigma\to0.
\end{equation}
Let us denote the coordinates of the point $\alpha^*$ by 
$\alpha_1^*=I_{nat}$, $\alpha^*_2=I_{int}$ and assume that they satisfy \eqref{II}.
Let $E_0^\sigma E_1^\sigma$ be the line segment of equilibrium states of system \eqref{r0model}
with the basic reproduction number given by 
\begin{equation}\label{pre'}
R_0(t)=R_0^{nat} - (R_0^{nat}-R_0^{int})\iint_{\Pi} q(\alpha;\sigma)\, \bigl({\mathcal R}_{\alpha}[r_{\alpha}^0]I\bigr)(t) \,d\alpha_1 d\alpha_2,\qquad t\ge 0,
\end{equation}
(cf.~\eqref{pre}).

\begin{proposition}\label{p2}
Let $E_0^\sigma E_1^\sigma$ be the line segment of equilibrium states of system \eqref{r0model}, \eqref{pre'}.
Let \eqref{sigma} hold where the coordinates of $\alpha^*=(I_{nat},I_{int})$ satisfy \eqref{II}.
Then, the following statements hold:
\begin{itemize}
\item[--] if the line segment $E_{int}E_{nat}$ defined by \eqref{EE} lies to the left of the line $I=I_{nat}$, then 
$E_0^\sigma, E_1^\sigma \to E_{nat}$ as $\sigma\to 0$;
\item[--] if the line segment $E_{int}E_{nat}$  lies to the right of the line $I=I_{int}$, then 
$E_0^\sigma, E_1^\sigma \to E_{int}$ as $\sigma\to 0$;
\item[--] if the line segment $E_{int}E_{nat}$ intersects the band $\mathcal I=\{(I,S): I_{nat}\le I\le I_{int}\}$, then 
$E_0^\sigma E_1^\sigma\to E_{int}E_{nat}\cap \mathcal I$
as $\sigma\to 0$.
\end{itemize}
\end{proposition}

\begin{figure}[ht!] 
    \centering
    \text{(a)}\includegraphics[width = 0.45\linewidth]{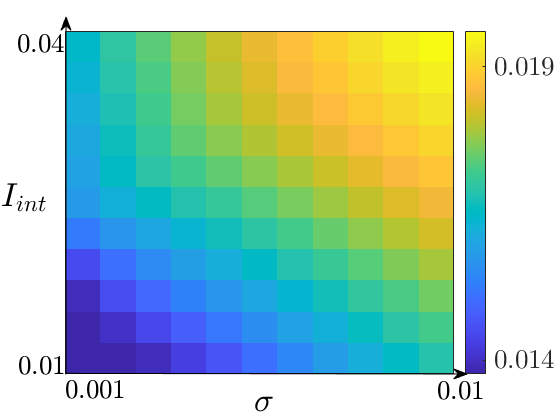}
    \text{(b)}\includegraphics[width = 0.45\linewidth]{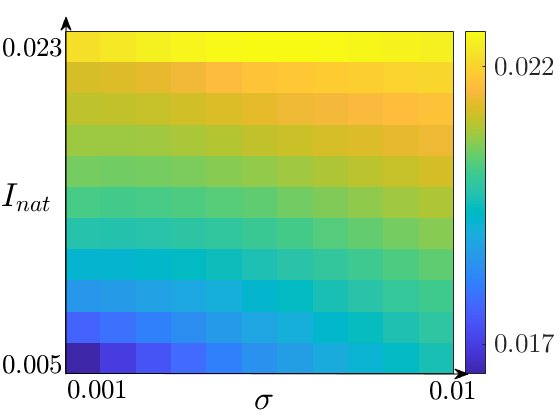}
    \caption{A trajectory of system \eqref{r0model}, \eqref{pre} with initial condition $I(0)=0$, $S(0)=1$
    converges to an endemic equilibrium $E_{\theta}=(I^*_{\theta},S^{*}_{\theta})$ which belongs to the line segment
    $E_0E_1$ of equilibrium states \eqref{endemic}. The color code presents the value $I^*_{\theta}$ of the infected population at the equilibrium as a function of parameters of the density function \eqref{density} of the Preisach operator \eqref{pre}. (a) Dependence of $I^*_{\theta}$ on $\sigma, I_{int}$ for fixed $I_{nat}=0.005$; (b) Dependence of $I^*_{\theta}$ on $\sigma, I_{int}$ for fixed $I_{int}=0.04$. 
    Other parameters are $\rho = 0.05$, $R_0^{nat} = 2.38$, and $R_0^{int} = 1.38$. The value of $I^*_{\theta}$ increases with $\sigma, I_{int}, I_{nat}$, i.e.~either delaying the intervention or revoking the intervention earlier or a higher degree of heterogeneity of the public response all result in a larger proportion of infected individuals at the endemic equilibrium state after the epidemic. 
    \label{prop2fig}}
\end{figure}

This alternative can be compared to Proposition \ref{p0} describing the endemic equilibrium states of switched system 
\eqref{r0model}, \eqref{relay1eq} depending on the relative position of points \eqref{EE} and the thresholds lines $I=I_{nat}, I_{int}$.
The proof is presented in Appendix 6.3.

\medskip\noindent{{\bf 3.5.} }
We are now interested in the global maximum value of the infected population density $I$
along the trajectory of the epidemic. 
This quantity is important for estimating the maximum burden on the healthcare system.
The trajectory of interest corresponds to the introduction
of a small number of infected individuals into an entirely susceptible population, i.e.~$I(0)\ll 1$ and 
$S(0)\approx 1$. 

\begin{figure}[ht!] 
    \centering
    \text{(a)}\includegraphics[width = 0.45\linewidth]{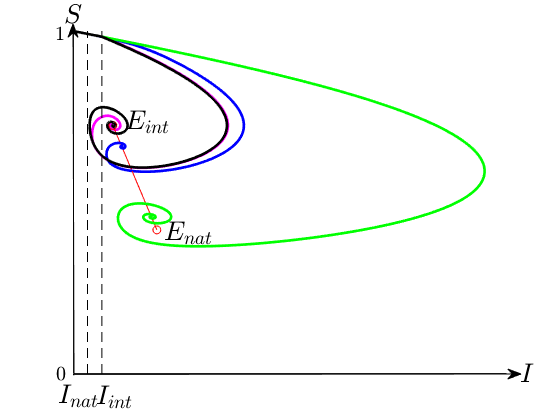}
    \text{(b)}\includegraphics[width = 0.45\linewidth]{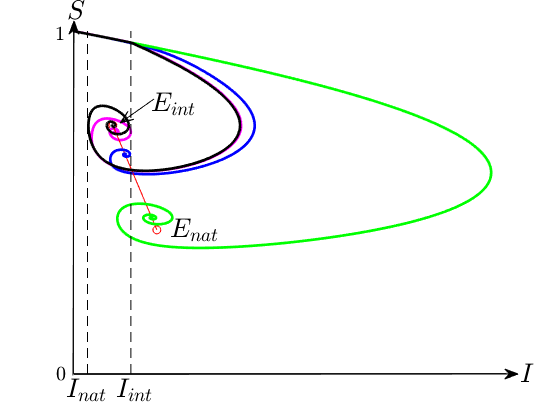}
    \text{(c)}\includegraphics[width = 0.45\linewidth]{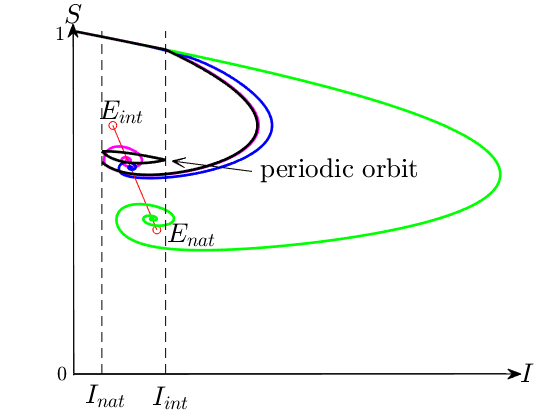}
   \text{(d)} \includegraphics[width = 0.45\linewidth]{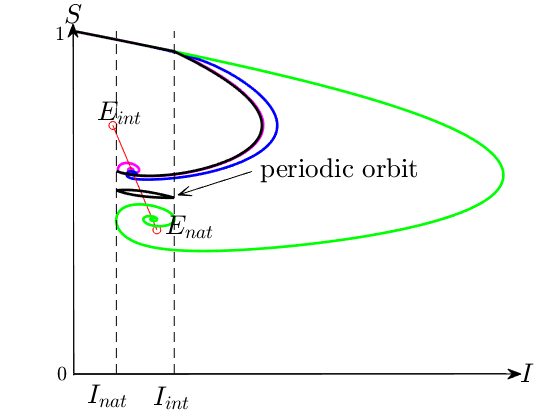}\hspace{1in}
   \text{(e)} \includegraphics[width = 0.45\linewidth]{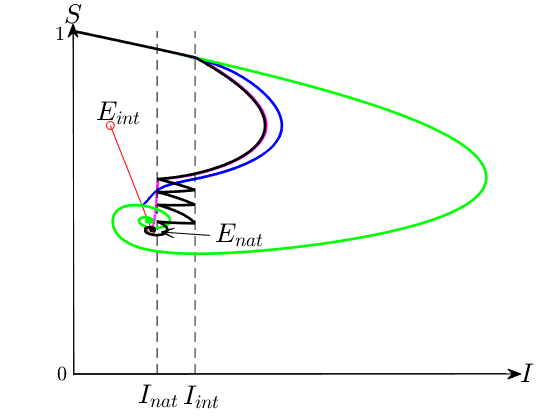}
   \text{(f)} \includegraphics[width = 0.45\linewidth]{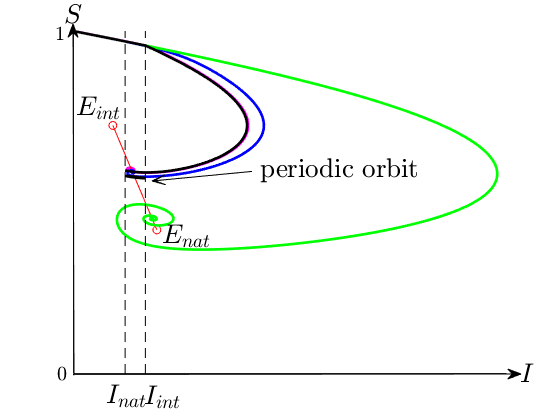}
    \caption{Convergence of a trajectory of system \eqref{r0model}, \eqref{pre} to an endemic equilibrium state $E_\theta$ for different values of the parameters $\sigma, I_{nat}, I_{int}$ of the density function \eqref{q}. 
    The point $E_\theta$ 
    belongs to the line segment $E_{int}E_{nat}$ 
    according to Propositions \ref{p1}, \ref{p2}.
    The pink, blue and green trajectories correspond the $\sigma$ values $0.001, 0.01$ and $0.1$ respectively. 
    The maximum of the infected population along the trajectory, $\max_{t\in\mathbb{R}_+} I(t)$, increases with 
    the heterogeneity parameter $\sigma$ for each pair $(I_{nat},I_{int})$. The black curves represent a trajectory of the ideally homogeneous model \eqref{r0model}, \eqref{relay1eq}. In panel c,d,f the homogeneous model converges to a periodic orbit.
      The initial values and parameters are $\rho = 0.05$, $R_0^{nat} = 2.38$, $R_0^{int} = 1.38$, $E_{int}= (0.013,0.73)$, $E_{nat}=(0.029,0.42)$ $I(0)=0.0001$, $S(0)=0.999$.
    \label{prop3fig}}
\end{figure}

Let us consider the ordinary differential system \eqref{r0model} where $R_0=\hat R^0(I)$ is a function of $I$ defined by
\begin{equation}\label{R00}
 \hat R^{0}(I) = R_0^{nat}-(R_0^{nat}-R_0^{int})\int_0^{I}\int_0^{\alpha_2}q(\alpha) \,d\alpha_1 d\alpha_2, \qquad R_0(t)=\hat R^0(I(t)).
\end{equation}
Denote by $\Gamma_0$ the unstable manifold of the saddle equilibrium $E_*=(0,1)$ of this system.
That is, $\Gamma_0=(I^0_\Gamma(t), S^0_\Gamma(t))$ ($t\in \mathbb{R}$) is the heteroclinic trajectory connecting the infection free equilibrium $E_*$
with the positive stable equilibrium $E_0=(I^*_0,S^*_0)$ of system \eqref{r0model}, \eqref{R00}. 
Assume that the equilibrium $E_0$ is a focus and denote by $M_0=(I_M^0,S_M^0)$ the first intersection point of the trajectory $\Gamma_0$ with the nullcline $\dot I=0$ of system \eqref{r0model}, \eqref{R00}. In other words,
\[
\max_{t\in\mathbb{R}} I_\Gamma^0(t) = I_M^0.
\]

\begin{proposition}\label{p3}
Suppose that the positive equilibrium $E_0$ of system \eqref{r0model}, \eqref{R00} is of focus type.
Then, trajectories $\mathcal C=(I(t),S(t))$ $(t\ge 0)$ of system \eqref{r0model}, \eqref{pre} with the Preisach operator
satisfy
\[
\sup_{t>0} I_\mathcal{C}(t)\to I_M^0 \quad \ \ \text{as} \quad \ \ I(0)\to 0, \ S(0)\to 1, \ I(0)+S(0)\le 1.
\]
This convergence is uniform with respect to an admissible initial state function of the Preisach operator \eqref{pre}.
\end{proposition}

\begin{figure}[ht!] 
    \centering
    \includegraphics[width = 0.7\linewidth]{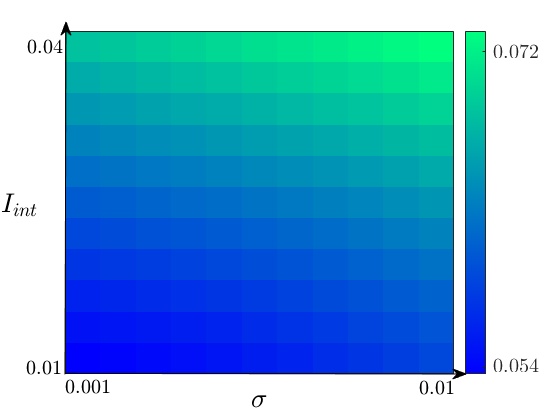}
    \caption{The color code represents the maximum of the infected population, $\max_{t\in\mathbb{R}_+} I(t)$,
along a trajectory of system \eqref{r0model}, \eqref{pre} for different values 
of the parameters $\sigma$ and $I_{int}$ of the density function \eqref{q} with fixed $I_{nat} = 0.005$. The maximum increases with $\sigma$ and $I_{int}$, i.e.~either delaying the intervention or a higher degree of heterogeneity of the public response result in a higher peak of infection during the endemic. The parameters of the system and the initial conditions of the trajectory are $\rho = 0.05$, $R_0^{nat} = 2.38$, $R_0^{int} = 1.38$, $E_{int}= (0.013,0.73)$, $E_{nat}=(0.029,0.42)$ $I(0)=0.0001$, $S(0)=0.999$ \label{prop4fig}}
\end{figure}

This proposition is used below to evaluate the maximum of $I$ numerically.
 The proof 
 is presented in Appendix 6.4.
 As shown in the proof, $\Gamma_0$ is a trajectory of both the ordinary differential 
 system \eqref{r0model}, \eqref{R00} and system \eqref{r0model}, \eqref{pre} with the Preisach operator; and, the endemic equilibrium $E_0$ of system \eqref{r0model}, \eqref{R00} coincides with the right end of the line segment $E_0E_1$ of endemic equilibrium states of system \eqref{r0model}, \eqref{pre}.
 
 Typically, the recovery rate $\gamma$ is much higher than the mortality rate $\mu$, hence $\rho\ll 1$. In this case, \eqref{focus} implies that the endemic equilibrium $E_0=(I_0^*,S_0^*)$ of system \eqref{r0model}, \eqref{R00} is a focus if the density $q=q(\alpha)$ is  close to a $\delta$-function. On the other hand, if $E_0$ is a node, then the heteroclinic trajectory $\Gamma_0$ converges to $E_0$
 without crossing the nullcline $\dot I=0$, therefore $\sup_{t\in\mathbb{R}} I_\Gamma^0(t)= I_0^*$. 
 One can show that in this case trajectories of system \eqref{r0model}, \eqref{pre} with the Preisach operator satisfy
 \[
 \sup_{t>0} I_\mathcal{C}(t)\to I_0^* \quad \ \ \text{as} \quad \ \ I(0)\to 0, \ S(0)\to 1, \ I(0)+S(0)\le 1.
 \]

\section{Discussion}
A number of observations can be made from numerical simulations of the heterogeneous model \eqref{r0model}, \eqref{pre}. 
As the density function $q$ in \eqref{pre}, we used the truncated Gaussian function
\begin{equation}\label{density}
q(\alpha)=Ae^{-\frac{(\alpha_1-I_{nat})^2+(\alpha_2-I_{int})^2}{2\sigma^2}}, \qquad \alpha=(\alpha_1,\alpha_2)\in 
\Pi=\{(\alpha_1,\alpha_2): 0\le \alpha_1<\alpha_2\le 1\},
\end{equation}
where $A$ is defined by the normalization condition \eqref{q}.
We observed the convergence to an endemic equilibrium state in all simulations. This is in contrast to switched system \eqref{r0model}, \eqref{relay1eq} where certain threshold pairs lead to a periodic orbit (see Figure \ref{relayfig}(abdf), \ref{colormap}).
According to Proposition \ref{p1},
a linear segment of endemic equilibrium states exists in the heterogeneous system. These endemic states differ by proportions of infected, susceptible and recovered individuals. The convergence upon  the segment of endemic equilibrium states depends on the degree of heterogeneity of the public response (measured by the parameter $\sigma$) and the threshold values at which the intervention is set to begin (the threshold $I_{int}$) and end (the threshold $I_{nat}$). In particular, Figure \ref{prop2fig}(a) shows that if the intervention starts later, i.e.~the threshold $I_{int}$ is higher, then the trajectory of the epidemic converges to an endemic state with higher proportions 
of infected and recovered individuals and a lower susceptible population. 
Revoking the intervention earlier, i.e.~increasing the threshold $I_{nat}$, 
has a similar effect (see Figure \ref{prop2fig}(b)). 
This is to be expected. More interestingly, the same figures demonstrate that more heterogeneity in the public response tends to steer the epidemic towards an endemic equilibrium state with larger infected and recovered populations.
However, there are exceptions to this general trend such as in Figure \ref{prop3fig}(e) where  the infected population at the endemic equilibrium state depends on the degree of heterogeneity of the public response, $\sigma$, in a non-monotonic fashion.
Further, we found that the maximum of the infected population along the epidemic trajectory demonstrates sensitive dependence on the heterogeneity parameter $\sigma$. Specifically, a higher heterogeneity of the public response leads to a higher peak of infection
(see Figures \ref{prop3fig}, \ref{prop4fig}).
These numerical results agree with, and complement, the statements of Propositions \ref{p1}\,--\,\ref{p3}.

The density function $q=q(\alpha)$ of the Preisach operator can be estimated from simultaneous observations of $I=I(t)$ and $R_0=R_0(t)$ using the Mayergoyz identification theorem (Mayergoyz 2003). If $I$ increases from zero to a value $I_1$, then drops back to zero, then increases to a value $I_2$, the drops back to zero again, etc., and if the local maximum values $\{I_1, I_2,\ldots, I_N\}$ of $I$ form an $\varepsilon$-net of the interval $0\le I\le 1$, then the identification theorem provides an $O(\varepsilon)$-approximation of the density function from measurements of $I$ and $v$. In the epidemiological context, this scenario corresponds to several waves of the epidemic. A number of practical identification algorithms can deal with measurement noise and limited amount of data, see e.g.~Hoffmann et al., 1989; Cirrincione et al., 2002; Rachinskii et al., 2016.
They include both non-parametric and parametric identification methods, where the latter assume a particular form of the density function such as in \eqref{density} or other (Appelbe et al., 2009; Krej\v{c}\'\i\ et al., 2011; Krej\v{c}\'\i\ et al., 2006; Brokate et al., 2011). 

\section{Conclusion}

We considered an SIR model where the transmission coefficient changes in response to
dynamics of the epidemic. We assumed that the adaptive response of an individual to the varying number of active cases can be modeled
by a two-state two-threshold hysteretic switch. In an ideally homogeneous population, the two switching thresholds of the transmission coefficient
can be imposed by 
the health authority which starts the intervention when the number of active cases exceeds a threshold $I_{int}$ and ends the intervention when the number of active cases drops below another threshold $I_{nat}$. 
In order to account for the possibility of a heterogeneous 
response among the susceptible individuals,
we allowed a distribution of switching thresholds and
modeled 
the aggregate response of the susceptible population by the Preisach operator.
The 
mean of the distribution represents the thresholds $I_{int}$ and $I_{nat}$ at which the health authority 
starts and ends the intervention;
the variance of the distribution, $\sigma^2$,
measures the degree of heterogeneity of the public response 
to the interventions. 
The resulting heterogeneous model 
is shown to have 
a continuum of endemic equilibrium states 
differing by the proportions of susceptible, infected and recovered populations.

Numerical simulations of the heterogeneous model provide an evidence that
a wider spread of thresholds of different population groups leads to a significant increase of the peak of infection during the epidemic. Further, a higher degree of heterogeneity of the public response tends to steer the epidemic trajectory to an endemic equilibrium state with higher proportions of the infected and recovered populations and a lower proportion of the susceptible population. In other words, a more uniform response of the public to transmission prevention measures helps ``flattening the curve" and can lead to smaller density of infection and lower $R_0$ 
when 
the endemic equilibrium state 
is reached
after the epidemic. These results suggest that 
intervention programs are more effective when accompanied by 
education campaigns 
which convince the public to comply with the intervention policies.
%
In particular, effective policy making  
should account for the culture and mindsets of 
the community when the intervention measures are decided upon,
and the importance of these measures should be conveyed to the public to ensure a more homogeneous response. 

The ideally homogeneous model predicts the lowest peak of infection.
On the other hand, some threshold pairs lead to the convergence of the epidemic trajectory to a periodic orbit
predicting  recurrent outbreaks of the epidemic in a homogeneous population when its response is hysteretic.
We observed that a slight degree of heterogeneity changes this scenario.
The heterogeneous response ensures the convergence to an equilibrium state after a higher infection peak.
These findings agree with the results of Kopfov\'a et al. (2021) where the heterogeneity of the response was shown to 
promote the global stability of the set of endemic equilibrium states in an SIR model with vaccination.
They 
are also associated with
a trade-off between achieving the herd immunity faster with a higher infection peak or in a controlled manner but slower.
Since expecting every member of the public to conform to the intervention policies exactly is 
unrealistic, some degree of heterogeneity with the associated effect on the epidemic dynamics can be assumed.
The optimal choice of thresholds of the intervention policy is an interesting problem which extends beyond the scope of this paper. In reality, this choice is constraint by many factors such as the capacity of the healthcare system, the cost
of the intervention measures to the economy and the willingness and readiness of the public to comply with intervention policies.




\section{Appendix}

{\bf 6.1.~Continuous Preisach model.}
Let us briefly recall 
a rigorous definition of the 
continuous Preisach operator \eqref{pre} (Krasnosel'skii et al.~1983). 
It involves a collection of non-ideal relays ${\mathcal R}_{\alpha}$, which respond to the same  continuous input $I=I(t)$ independently according to formula \eqref{relay'}. The relays contributing to the system have different pairs of thresholds $\alpha=(\alpha_1,\alpha_2)\in \Pi$,
where the subset $\Pi$ of the half-plane $\{\alpha=(\alpha_1,\alpha_2): \alpha_1<\alpha_2\}$ 
is assumed to be measurable and bounded;
the $\alpha$-plane is called the Preisach plane.
The output 
of the continuous Preisach model is the scalar-valued function $R_0=R_0(t): \mathbb{R}_+\to \mathbb{R}$
defined by \eqref{pre}, 
where 
$q=q(\alpha):\Pi\to\mathbb{R}_+$  is a positive bounded measurable function (measure density) representing the weights of the relays; and, $r_{\alpha}^0$
is the initial state of the relay ${\mathcal R}_{\alpha}$ for any given $\alpha\in\Pi$.  
The function $r^0=r_{\alpha}^0:\Pi\to \{0,1\}$ of the variable $\alpha=(\alpha_1,\alpha_2)$
is referred to as the {\em initial state} function of the Preisach operator. It
is assumed to be measurable and satisfy the constraints \eqref{v1}, \eqref{v2}, in which case the initial state-input pair is called {\em compatible}.
These requirements ensure that the integral in \eqref{pre} is well-defined for each $t\ge 0$ and, furthermore, the output $R_0=R_0(t)$ of the Preisach model is a continuous function of time.
The function \eqref{relay'} with a fixed $t\ge0$ and varying $\alpha\in \Pi$ is interpreted as the {\em state function} of the Preisach model at the moment $t$ as it describes the states of all the relays at this moment; this state function $r(t)=r_\alpha(t):\Pi\to\{0,1\}$ is an element of the space $L_1(\Pi;\mathbb{R})$ for each $t\ge0$.

For brevity, let us denote the input-to-output operator of the Preisach model defined by \eqref{pre} by
\begin{equation}\label{P}
R_0(t)=({\mathcal P}[r^0]I)(t),\qquad t\geq 0,
\end{equation}
where both the input $I:\mathbb{R}_+\to \mathbb{R}$ and the initial state function $r^0=r^0_{\alpha}$ (which is compatible with the input)
are the arguments; the value of this operator is the output
$R_0:\mathbb{R}_+\to \mathbb{R}$. 
%
%
%
%
An important property of the Preisach operator \eqref{pre} is that it is 
Lipschitz continuous if $q:\Pi\to\mathbb{R}_+$ is bounded (Krasnosel'skii et al. 1983). More precisely, the relations 
\[
R_0^k(t)=({\mathcal P}[r^{0,k}]I^k)(t),\qquad t\geq 0, {\qquad k=1,2,}
\]
and $0\le I^1(t), I^2(t)\le 1$ ($t\ge 0$) imply
\begin{equation}\label{LipP*}
\|R_0^1-R_0^2\|_{C([0,\tau];\mathbb{R})} \le  K \Big( \|r^{0,1}-r^{0,2}\|_{L_1(\Pi;\mathbb{R})} + \|I^1-I^2\|_{C([0,\tau];\mathbb{R})}\Big)
\end{equation}
for any $\tau\geq 0$ with the Lipschitz constant
\begin{equation}\label{K}
K := \max_{0\le \alpha_1<\alpha_2\le 1} q(\alpha).
\end{equation}

	Let us denote by $\mathfrak U$ the set of all triplets $(I_0,S_0,r^0)$, where $(I_0,S_0)\in \mathfrak D=\{(I,S): 0\le I, S; I+S\le 1\}$ 
	and the initial state function $r^0=r^0_\alpha$ of the Preisach operator is compatible with $I_0$. 
	The global Lipschitz estimate \eqref{LipP*} ensures (for example, using the Picard-Lindel\"of type of argument) that for a
	given $(I_0,S_0,r^0)\in \mathfrak U$, 
	system \eqref{r0model} with the Preisach operator \eqref{pre} has a unique local solution
	with the initial data $I(0)=I_0, S(0)=S_0$ and the initial state function $r^0$ (see, for example, the survey (Leonov 2017)). Further, the positive invariance of $\mathfrak D$ implies that each solution is extendable to the whole semi-axis $t\geq 0$. These solutions induce a continuous {semi-flow} in the set $\mathfrak U$, which is considered to be the phase space 
	of system \eqref{r0model}, \eqref{pre} and is endowed with a metric by the natural embedding into the space $\mathbb{R}^2\times L_1(\Pi;\mathbb{R})$. This construction leads to the standard definition of local and global stability including stability of equilibrium states and periodic solutions. In particular, an equilibrium is a triplet $(I_0,S_0,r^0)\in \mathfrak{U}$ and a periodic solution is a periodic function $(I(\cdot),S(\cdot),r(\cdot)): \mathbb{R}_+\to \mathfrak U$ where the last component, viewed as a function $r:\mathbb{R}_+\times \Pi\to \{0,1\}$ of two variables $t\in\mathbb{R}_+$ and $\alpha\in\Pi$, is given by \eqref{re}.
The basic reproduction number \eqref{pre} at an equilibrium is constant, while for a periodic solution the basic reproduction number changes periodically with the period of $I$ and $S$.

\medskip
\noindent
{\bf 6.2.~Proof of Proposition 2.}
%
Consider the strictly decreasing functions \eqref{R00} and
\begin{equation}\label{R01}
    \hat R^{1}(I) = R_0^{int}+(R_0^{nat}-R_0^{int})\int_{I}^{1}\int_I^{\alpha_2}q(\alpha)\,d\alpha_1 d\alpha_2.
\end{equation}
Due to the assumption $q(\alpha)>0$ and the normalization condition \eqref{q}, 
\begin{equation}\label{bla}
\hat R^0(I)> \hat R^1(I)\quad \text{for}\quad I\in(0,1); \quad
\hat R^{0}(0)=\hat R^{1}(0)=R_0^{nat};\quad \hat R^{0}(1)=\hat R^{1}(1)=R_0^{int}.
\end{equation}
For every $\theta \in [0,1]$ set
\begin{equation}\label{r0theta1}
    \hat R^{\theta}(I) = \theta \hat R^{1}(I) + (1-\theta) \hat R^{0}(I).
\end{equation}
The constraint \eqref{v2} implies that the basic reproduction number $R_0$
at an endemic equilibrium state $(I^*,S^*)$ equals  $\hat R^\theta(I^*)$ for some $\theta\in [0,1]$.
On the other hand, $I^*$ and $S^*$ are related to $R_0$
by formulas \eqref{equil}. Hence, $(I^*,S^*)$ is an endemic equilibrium state iff there is a $\theta\in[0,1]$ such that
\begin{equation}\label{nonl}
    I^*=\left(1-\frac{1}{\hat R^\theta(I^*)}\right)\rho
\end{equation}
and 
\begin{equation}\label{nonls}
    S^*=\frac{1}{\hat R^\theta(I^*)}.
\end{equation}

By definition, the positive function \eqref{r0theta1} is continuous in $I, \theta$, strictly decreases in $I$ and strictly decreases in $\theta$ for $I\in (0,1)$. Therefore, the function
\begin{equation}\label{f}
f_\theta(I)=I-\left(1-\frac{1}{\hat R^\theta(I)}\right)\rho
\end{equation}
strictly increases in $I$ and $\theta$. Since
\[
f_\theta(0)=-\left(1-\frac{1}{R_0^{nat}}\right)\rho<0,\qquad f_\theta(I)=1-\left(1-\frac{1}{R_0^{int}}\right)\rho>0, \qquad \theta\in[0,1],
\]
equation \eqref{nonl} has a unique solution $I^*=I^*_\theta\in (0,1)$ for each $\theta\in[0,1]$.
By the Implicit Function Theorem, $I^*_\theta$ depends continuously on, and strictly decreases in, $\theta$.
Therefore, formulas \eqref{nonl}, \eqref{nonls} imply that the endemic equilibrium states $(I^*_\theta,S^*_\theta)$
form a line segment $E_0E_1$ on the $(I,S)$-plane with
$E_0=(I^*_0,S^*_0)$ and $E_1=(I^*_1,S^*_1)$. 
Finally, from
\eqref{bla}, \eqref{r0theta1} it follows that $R_0^{int}\le \hat R^\theta(I)\le R_0^{nat}$, hence
from \eqref{nonl}, \eqref{nonls} and \eqref{EE} we conclude that
$E_0E_1\subset E_{int}E_{nat}$, which completes the proof.

%

\medskip
\noindent
{\bf 6.3.~Proof of Proposition 3.}
Condition \eqref{sigma} implies that given any $\varepsilon>0$ there is a $\sigma_0=\sigma_0(\varepsilon)$ such that for all $\sigma\in(0,\sigma_0)$ the following relations hold:
\begin{eqnarray}
\hat R^0(I)>R_0^{nat}-\varepsilon & \text{if} & I\le I_{int}-\varepsilon, \label{c1}\\
\hat R^0(I)<R_0^{int}+\varepsilon & \text{if} & I\ge I_{int}+\varepsilon, \label{c2}\\
\hat R^1(I)>R_0^{nat}-\varepsilon & \text{if} & I\le I_{nat}-\varepsilon, \label{c3}\\
\hat R^1(I)<R_0^{int}+\varepsilon & \text{if} & I\ge I_{nat}+\varepsilon. \label{c4}
 \end{eqnarray}
 Without loss of generality, we assume that
 \begin{equation}\label{bc}
 0<\varepsilon<R_0^{nat}-1.
 \end{equation}
 
Let us consider different cases of positioning he threshold $I_{int}$.
 
 \smallskip
 {\bf Case 1:} $I_{int}\le \rho(1-1/R_0^{int})
 $. In this case, set
 \[
 I_-=\rho(1-1/R_0^{int})
 , \qquad I_+=\rho(1-1/R_0^{int})+\varepsilon
 \]
 and notice that
$
 I_+ 
 \ge I_{int}+\varepsilon
$
implies $\hat R^0(I_+)<R_0^{int}+\varepsilon$ due to \eqref{c2}. Therefore, the function \eqref{f} satisfies
\begin{equation}\label{a11}
f_0(I_+)=I_+-\rho(1-1/\hat R^0(I_+))=\varepsilon-\frac{\rho}{R_0^{int}}+\frac{\rho}{\hat R^0(I_+)}>
\varepsilon-\frac{\rho}{R_0^{int}}+\frac{\rho}{R_0^{int}+\varepsilon}>0
\end{equation}
because $\rho<1<R_0^{int}$.
On the other hand, since $R_0^{int}\le \hat R^0(I)$ for all $I\in[0,1]$,
\begin{equation}\label{a12}
f_0(I_-)=\rho(1-1/R_0^{int})-\rho(1-1/\hat R^0(I_-))\le 0.
\end{equation}
Hence, a unique root $I_0^*$ of the equation $f_0(I)=0$ satisfies $I_-\le I_0^*<I_+$, i.e.
\begin{equation}\label{b1}
\rho(1-1/R_0^{int})\le I_0^*<\rho(1-1/R_0^{int})+\varepsilon.
\end{equation}

\medskip
 {\bf Case 2:} $\rho(1-1/R_0^{int})< I_{int}< \rho(1-1/R_0^{nat}) 
 $. Set $I_-=I_{int}-\varepsilon$, $I_+=I_{int}+\varepsilon$. Then,
 \[
 f_0(I_+)=I_{int}+\varepsilon-\rho(1-1/\hat R^0(I_+))>\varepsilon-\frac{\rho}{R_0^{int}}+\frac{\rho}{\hat R^0(I_+)}>0
 \]
 as in \eqref{a11}. On the other hand, 
 \[
 f_0(I_-)=I_{int}-\varepsilon-\rho(1-1/\hat R^0(I_-))<-\varepsilon -\frac{\rho}{R_0^{nat}}+\frac{\rho}{\hat R^0(I_-)}
 \]
 and due to \eqref{c1},
 \begin{equation}\label{cd'}
 f_0(I_-)
 <-\varepsilon -\frac{\rho}{R_0^{nat}}+\frac{\rho}{R_0^{nat}-\varepsilon}<0,
 \end{equation}
 where the last inequality follows from \eqref{bc}.
 Therefore, the root $I_0^*$ of $f_0(I)=0$ satisfies
 \begin{equation}\label{b2}
|I_0^*-I_{int}|<\varepsilon
.
\end{equation}

\medskip
 {\bf Case 3:} $I_{int}\ge \rho(1-1/R_0^{nat})
 $. In this case, set
 \[
 I_-=\rho(1-1/R_0^{nat})-\varepsilon
 , \qquad I_+=\rho(1-1/R_0^{nat})
 \]
 and notice that
$
 I_- 
 \le I_{int}-\varepsilon
$
implies $\hat R^0(I_-)>R_0^{nat}-\varepsilon$ due to \eqref{c1}. Therefore,
\[
f_0(I_-)=\rho(1-1/R_0^{nat})-\varepsilon-\rho(1-1/\hat R^0(I_-))<-\varepsilon -\frac{\rho}{R_0^{nat}}+\frac{\rho}{R_0^{nat}-\varepsilon}<0
\]
as in \eqref{cd'}. On the other hand, since $\hat R_0(I)\le R_0^{nat}$ for all $I\in[0,1]$,
\[
f_0(I_+)=\rho(1-1/R_0^{nat})-\rho(1-1/\hat R^0(I_+))\ge 0,
\]
hence $I_-<I_0^*\le I_+$, which is equivalent to
 \begin{equation}\label{b3}
\rho(1-1/R_0^{nat})-\varepsilon<I_0^*\le \rho(1-1/R_0^{nat}).
\end{equation}

\medskip
Since $\varepsilon>0$ is arbitrarily small, relations \eqref{b1}, \eqref{b2}, \eqref{b3} obtained in Cases 1\,--\,3, respectively, imply that the component $I_0^*=I^*_0(\sigma)$ of the equilibrium state $E_0^\sigma=(I^*_0(\sigma),S^*_0(\sigma))$ satisfies
\begin{equation}\label{I*0}
I^*_0(\sigma) \to 
\begin{cases}
\rho(1-1/R_0^{int}) & \text{if} \ \ \ I_{int}\le \rho(1-1/R_0^{int}),\\
I_{int} & \text{if} \ \ \ \rho(1-1/R_0^{int})<I_{int}< \rho(1-1/R_0^{nat}),\\
\rho(1-1/R_0^{nat}) & \text{if} \ \ \ I_{int}\le \rho(1-1/R_0^{nat})
\end{cases}
\end{equation}
as $\sigma\to 0$. The same argument shows that relations \eqref{c3}, \eqref{c4} imply a similar relationship
\begin{equation}\label{I*1}
I^*_1(\sigma) \to 
\begin{cases}
\rho(1-1/R_0^{int}) & \text{if} \ \ \ I_{nat}\le \rho(1-1/R_0^{int}),\\
I_{nat} & \text{if} \ \ \ \rho(1-1/R_0^{int})<I_{nat}< \rho(1-1/R_0^{nat}),\\
\rho(1-1/R_0^{nat}) & \text{if} \ \ \ I_{nat}\le \rho(1-1/R_0^{nat})
\end{cases}
\end{equation}
as $\sigma\to 0$ for the equilibrium state $E_1^\sigma=(I^*_1(\sigma),S^*_1(\sigma))$.
Relations \eqref{I*0}, \eqref{I*1} are equivalent to the conclusion of the proposition.

\medskip
\noindent
{\bf 6.4.~Proof of Proposition 4.} 
Consider the ordinary differential system \eqref{r0model} where $R_0=\hat R^0(I)$:
\begin{equation}\label{r0model'}
\begin{aligned}
    \dot I &= \hat R^0(I) SI - I, \\
    \dot S &= - \hat R^0(I) SI - \rho S + \rho.
\end{aligned}
\end{equation}
The positive equilibrium $E_0=(I_0^*,S_0^*)$ of this system coincides with the right end of the line segment of equilibrium states of system \eqref{r0model}, \eqref{pre}.

\begin{figure}[ht!]
    \centering
    \includegraphics[width = 0.6\linewidth]{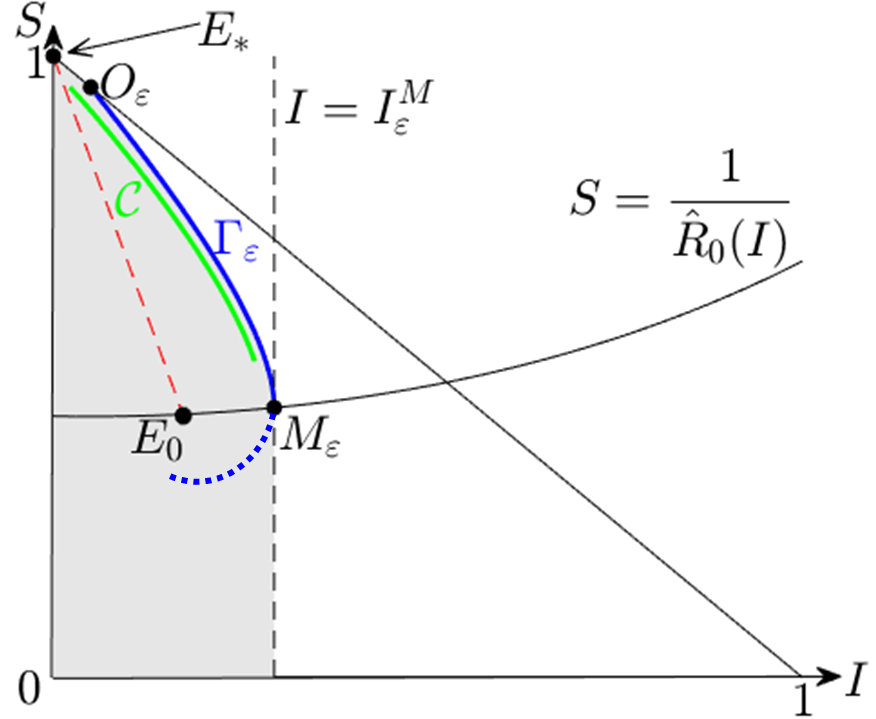}
    \caption{Points $E_* = (0,1)$, $O_{\varepsilon} = (I_0^{\varepsilon},S_0^{\varepsilon})$, $M_{\varepsilon} = (I_M^{\varepsilon}, S_M^{\varepsilon})$, $E_0 = (I_0^*,S_0^*)$, the nullcline $S = {1}/{\hat R_0 (I)}$ and the trajectories $\Gamma_{\varepsilon}$ and $\mathcal{C}$. The shaded domain $Q_\varepsilon$ is invariant for system \eqref{r0model}, \eqref{pre}.}
    \label{prop4}
\end{figure}

Take an $\varepsilon>0$ and denote $O_\varepsilon=(I_O^\varepsilon,S_O^\varepsilon)=(\varepsilon,1-\varepsilon)$.
Denote by $\hat\Gamma_\varepsilon$ the trajectory of \eqref{r0model'} starting at the point $O_\varepsilon$
and continued in forward time. 
Let $\varepsilon>0$ be sufficiently small to ensure that the point $O_\varepsilon$ lies above the nullcline
$\dot I=0$ of system \eqref{r0model'} defined by $S=1/\hat R^0(I)$.
By assumption, the positive equilibrium $E_0$ of \eqref{r0model'} is a focus.
Therefore, the trajectory $\hat \Gamma_\varepsilon$ intersects the nullcline $S=1/\hat R^0(I)$.
Denote by $M_\varepsilon=(I_M^\varepsilon,S_M^\varepsilon)$ the first point of intersection of $\hat \Gamma_\varepsilon$ with the nullcline $S=1/\hat R^0(I)$ 
and by 
$\Gamma_\varepsilon$ the arc $O_\varepsilon M_\varepsilon$ of the curve $\hat \Gamma_\varepsilon$,
see Figure \ref{prop4}.
Since $\dot I>0$ above the nullcline $S=1/\hat R^0(I)$, the curve $\Gamma_\varepsilon$ is the graph 
of a continuous function $\gamma_\varepsilon: [I_O^\varepsilon,M_O^\varepsilon]\to\mathbb{R}$:
\[
\Gamma_\varepsilon=\{(I,S): S=\gamma_\varepsilon(I), \ I\in [I_O^\varepsilon,M_O^\varepsilon]\}.
\]

The main step of the following proof is to show that a trajectory $\mathcal C=(I_\mathcal{C}(t),S_\mathcal{C}(t))$ ($t\ge0$) of system \eqref{r0model}, \eqref{pre}, which starts 
below the curve $E_* O_\varepsilon\cup \Gamma_\varepsilon$ in the band $\{(I,S): 0<I<I_M^\varepsilon\}$
with an arbitrary admissible initial state function $r_\alpha^0$
of the Preisach operator, does not cross the line $I=I_M^\varepsilon$ when continued in forward time. To this end,
let us show that if $\mathcal C \cap \{(I,S): I=I_M^\varepsilon\}=\varnothing$, then the following implication holds:
\begin{equation}\label{imp}
I_\mathcal{C}(t)< I_\mathcal{C}(\hat t)=I_M^\varepsilon \ \ \text{for all} \ \ 0\le t < \hat t \ \ \Rightarrow \ \
(I_\mathcal{C}(\hat t),S_\mathcal{C}(\hat t))=M_\varepsilon.
\end{equation}

Assume that there is a $\hat t>0$ such that 
$I_\mathcal{C}(t)< I_\mathcal{C}(\hat t)=I_M^\varepsilon$ for $t \in [0, \hat t)$. 
If we assume that $S_\mathcal{C}(\hat t)<S_M^\varepsilon$, then $S_\mathcal{C}(t)<1/\hat R^0(I_\mathcal{C}(t))$ on a sufficiently small time interval $[\hat t-\delta,\hat t]$ 
 and hence
\begin{equation}\label{bll}
0<I_\mathcal{C}(\hat t) -I_\mathcal{C}(\hat t-\delta)=\int_{\hat t-\delta}^{\hat t}\!\! \dot I_\mathcal{C}(t)dt
=\int_{\hat t-\delta}^{\hat t} \!\!I_\mathcal{C}(t)\bigl( R_0(t) S_\mathcal{C}(t)-1\bigr)dt
< \int_{\hat t-\delta}^{\hat t} \!\!I_\mathcal{C}(t)\!\!\left( \frac{R_0^\mathcal{C}(t)}{\hat R^0(I_\mathcal{C}(t))}-1\right)\!dt,
\end{equation}
where the basic reproduction number $R_0^\mathcal{C}(t)=({\mathcal P}[r^0]I_\mathcal{C})(t)$ is given by \eqref{pre} with $I=I_\mathcal{C}(t)$.
But equations \eqref{compatibility} and \eqref{pre} imply that
\begin{equation}\label{bl}
R_0(t)\le \hat R^0(I(t)); \qquad \quad R_0(t)= \hat R^0(I(t)) \ \ \Leftrightarrow \ \ r_\alpha(t)=\begin{cases} 1 \ \ \text{if} \ \ \alpha_2<I(t),\\
0 \ \ \text{if} \ \  \alpha_2>I(t),
\end{cases}
\end{equation}
\begin{equation}\label{rrr}
   R_0^{int} \le R_0(t)\le R_0^{nat} 
\end{equation}
at any time $t\ge0$ for any input and initial state function of the Preisach operator. In particular,
$R_0^\mathcal{C}(t)\le \hat R^0(I_\mathcal{C}(t))$ at all times, hence 
\begin{equation}\label{bnn}
I_\mathcal{C}(t)\left( \frac{R_0^\mathcal{C}(t)}{\hat R^0(I_\mathcal{C}(t))}-1\right) 
\le 0,
\end{equation}
which is in contradiction with \eqref{bll}. This contradiction shows that $S_\mathcal{C}(\hat t)\ge S_M^\varepsilon$.
Since the initial point of the trajectory $\mathcal{C}$ lies below the curve $E_*O_\varepsilon \cup\Gamma_\varepsilon $ and 
the curve $\mathcal C$ cannot cross the line segment $E_*O_\varepsilon$
due the positive invariance of the
domain $0\le I,S,\ I+S\le 1$, from $S_\mathcal{C}(\hat t)\ge S_M^\varepsilon$ we conclude that 
$\mathcal{C}\cap \Gamma_\varepsilon\ne \varnothing$ and there is a $t'\in (0,\hat t]$ such that
\begin{equation}\label{intersect}
I_{\mathcal{C}}(t)< I_M^\varepsilon, 
\ \ (I_{\mathcal{C}}(t), S_{\mathcal{C}}(t))\not\in \Gamma_\varepsilon 
\ \ \ {\rm for} \ \ \ {0\le t<t'\le \hat t}; 
\qquad (I_{\mathcal{C}}(t'),S_{\mathcal{C}}(t'))\in \Gamma_\varepsilon.
\end{equation}

%
Next, let us
consider the determinant
\begin{equation}\label{Delta}
\Delta=\left|\begin{array}{ccc}
\hat R^0(I) SI-I &&  R_0 SI-I\\
-\hat R^0(I)SI+(1-S)\rho && - R_0SI+(1-S)\rho
\end{array}
\right|=SI\bigl(\hat R^0(I)-R_0\bigr)\bigl( (1-S)\rho-I\bigr),
\end{equation}
where the columns are the vector fields of systems \eqref{r0model'} and \eqref{r0model}, respectively.
Notice that the segment $E_* E_0$, which belongs to the straight line $(1-S)\rho-I=0$,
lies above the curve $S=1/\hat R^0(I)$ (because $S$ increases along this curve) and meets it at the point $E_0$,
see Figure \ref{prop4}.
On this segment, with the exception of the point $E_0=(I_0^*,S_0^*)$, the vector field of system \eqref{r0model'} satisfies
\[
\frac{dS}{dI}=\frac{-\hat R^0(I)SI+(1-S)\rho}{\hat R^0(I) SI-I}=-1>-\frac1{\rho}.
\]
Since the slope of the line segment $E_* E_0$ is $-1/\rho$ and the trajectory $\Gamma_\varepsilon$ starts at the point $O_\varepsilon$ above this segment, it follows that $\Gamma_\varepsilon$ lies strictly above $E_* E_0$  in the vertical band $I_O^\varepsilon\le I\le I_0^*$. Further, the line $(1-S)\rho-I=0$ lies strictly below the curve $S=1/\hat R^0(I)$
for $I>I_0^*$, and $\Gamma_\varepsilon$ lies above this curve,
hence $\Gamma_\varepsilon$ lies strictly above the straight line $(1-S)\rho-I=0$ for all $I\in [I_O^\varepsilon,M_O^\varepsilon]$.
Therefore, \eqref{Delta} implies 
\begin{equation}\label{Delta'}
{\rm sign}\, \Delta ={\rm sign}\, \bigl(R_0-\hat R^0(I)\bigr)
\end{equation}
on the curve $\Gamma_\varepsilon$, and from \eqref{bl} it follows that
$\Delta\le 0$ on $\Gamma_\varepsilon$. On the other hand, 
the vector fields of systems \eqref{r0model'} and \eqref{r0model} should satisfy the opposite inequality $\Delta\ge0$ at the first intersection point $(I_\mathcal{C}(t'),S_\mathcal{C}(t'))$ of the trajectories $\mathcal{C}$ and $\Gamma_\varepsilon$ (cf.~\eqref{intersect}) because the initial point of the trajectory $\mathcal{C}$ lies to the left of the curve $E_*O_\varepsilon \cup\Gamma_\varepsilon$. Hence, 
$\Delta=0$ at the first intersection point, which due to \eqref{bl}, \eqref{Delta'} is equivalent to the relations
\begin{equation}\label{bc1}
R_0^\mathcal{C}(t')=\hat R^0(I_\mathcal{C}(t')), \qquad r_\alpha^\mathcal{C}(t')=\begin{cases} 1 \ \ \text{if} \ \ \alpha_2<I_\mathcal{C}(t')\\
0 \ \ \text{if} \ \  \alpha_2>I_\mathcal{C}(t')\end{cases}
\end{equation}
for the trajectory $\mathcal{C}$, where the basic reproduction number $R_0^\mathcal{C}(t)$ and the state function $r_\alpha^\mathcal{C}(t)$
are given by equations \eqref{pre} and \eqref{re} with $I=I_\mathcal{C}(t)$.

The definitions \eqref{relay'} and \eqref{pre} of the relay and the Preisach operator imply that $\Gamma_\varepsilon$
is a trajectory of a solution $(I_\Gamma^\varepsilon(t),S_\Gamma^\varepsilon(t))$ ($0\le t\le t_M^\varepsilon$) of system \eqref{r0model}, \eqref{pre}
and that the basic reproduction number and the state function for this solution are defined by
\begin{equation}\label{bc2}
R_0^\Gamma(t)=\hat R^0(I^\varepsilon_\Gamma(t)), \qquad r_\alpha^\Gamma(t)=\begin{cases} 1 \ \ \text{if} \ \ \alpha_2<I^\varepsilon_{\Gamma}(t)\\
0 \ \ \text{if} \ \  \alpha_2>I^\varepsilon_{\Gamma}(t)\end{cases}
\end{equation}
on the interval $0\le t\le t_M^\varepsilon$. Due to the forward uniqueness property of solutions of system \eqref{r0model},
\eqref{pre}, from \eqref{bc1}, \eqref{bc2} it follows that the trajectories $\mathcal{C}$ and $\Gamma_{\varepsilon}$
of this system coincide after the moment $t'$ when they merge, which completes the proof of the implication \eqref{imp}.

We conclude that if a trajectory $\mathcal{C}$ starting below the curve  $E_* O_\varepsilon\cup \Gamma_\varepsilon$  in the band $0<I<I_M^\varepsilon$ ever reaches the vertical line $I=I_M^\varepsilon$, then the intersection occurs at the point $M_\varepsilon$ of the curve $S=1/\hat R^0(I)$ where $\dot I_\mathcal{C}=0$, $\dot S_\mathcal{C}<0$. Hence, $\mathcal C$ crosses
the line $S=1/\hat R^0(I)$ at the point $M_\varepsilon$ vertically downwards to the domain $S<1/\hat R^0(I)$ where
\[
\dot I_\mathcal{C}(t)= I_\mathcal{C}(t) \bigl( R_0^\mathcal{C}(t) S_\mathcal{C}(t)-1\bigr)<I_\mathcal{C}(t)\left( \frac{R_0^\mathcal{C}(t)}{\hat R^0(I_\mathcal{C}(t))}-1\right)\le 0
\]
due to \eqref{bnn}. In other words, the point $\hat t$ defined by 
\eqref{imp}
is an isolated point of local 
maximum
for $I_\mathcal{C}(t)$. Thus, the domain $\mathcal Q_\varepsilon$ of the positive quadrant $I,S> 0$ bounded form above by the curve $E_* O_\varepsilon\cup \Gamma_\varepsilon$ and from the right by the vertical line $I=I_M^\varepsilon$ (see Figure \ref{new}) is invariant for  system \eqref{r0model}, \eqref{pre} with any admissible initial state function of the Preisach operator:
\[
(I_\mathcal{C}(0), S_\mathcal{C}(0))\in \mathcal{Q}_\varepsilon \quad \Rightarrow \quad (I_\mathcal{C}(t), S_\mathcal{C}(t))\in \mathcal{Q}_\varepsilon, \ \ t\ge0.
\]
In particular, any trajectory $\mathcal C$ of \eqref{r0model}, \eqref{pre} starting
in the $\varepsilon$-neighborhood of the point $E_*=(0,1)$
in the positively invariant domain $0\le I,S,\ I+S\le 1$ satisfies
\begin{equation}\label{III}
I_\mathcal{C}(t)\le I_M^\varepsilon, \ \ t\ge0.
\end{equation}

It remains to show that 
$\sup_{t>0} I_\mathcal{C}(t)\to I_M^0$ uniformly with respect to the initial condition
$(I_\mathcal{C}(0), S_\mathcal{C}(0))\in U_\varepsilon=\{(I,S): 0<I\le 1-S<\varepsilon\}$ and an admissible initial state function of the Preisach operator as $\varepsilon\to 0$. To this end, take $s_0, i_0, \varepsilon_0>0$ such that
\begin{equation}\label{bb1}
s_0>\frac1{R_0^{int}},
\end{equation}
\begin{equation}\label{bb2}
i_0<\left(1-\frac1{R_0^{int}s_0}\right)(1-s_0-\varepsilon_0),
\end{equation}
\begin{equation}\label{bb3}
i_0<\frac{(R_0^{int}s_0-1)\rho}{R_0^{nat}}.
\end{equation}
Take an $\varepsilon\in (0,\varepsilon_0)$ and consider a trajectory $\mathcal C$ of system \eqref{r0model}, \eqref{pre} satisfying $(I_\mathcal{C}(0), S_\mathcal{C}(0))\in U_\varepsilon$. 
Notice that at the initial moment $S_\mathcal{C}(0)>s_0$ and 
due to \eqref{rrr}, \eqref{bb1},
\begin{equation}\label{hh}
S_\mathcal{C}(t)>s_0, \qquad R_0^\mathcal{C}(t) S_\mathcal{C}(t)>R_0^{int}s_0>1 
\end{equation}
on some time interval $0\le t<\delta$. Further, as long as \eqref{hh} holds, also
$\dot I_\mathcal{C}(t)=I_\mathcal{C}(t)\bigl(R_0^\mathcal{C}(t) S_\mathcal{C}(t)-1\bigr)>0$ and
\[
\frac{\dot S_\mathcal{C}(t)}{\dot I_\mathcal{C}(t)} =-1+\frac{(1-S_\mathcal{C}(t))\rho-I_\mathcal{C}(t)}{I_\mathcal{C}(t)(R_0^\mathcal{C}(t) S_\mathcal{C}(t)-1)}\ge 
-1-\frac{1}{R_0^\mathcal{C}(t) S_\mathcal{C}(t)-1}>
-1-\frac1{R_0^{int}s_0-1},
\]
which implies
\[
S_\mathcal{C}(t)>S_\mathcal{C}(0)-\left(1+\frac1{R_0^{int}s_0-1}\right)(I_\mathcal{C}(t)-I_\mathcal{C}(0))>1-\varepsilon_0-\left(1+\frac1{R_0^{int}s_0-1}\right)I_\mathcal{C}(t).
\]
Due to this estimate and \eqref{bb2}, we conclude that the trajectory $\mathcal{C}$ reaches the line $I=i_0$ before it reaches the line $S=s_0$, i.e.~there is a $t'>0$ such that
\[
S_\mathcal{C}(t')=i_0, \qquad \min_{0\le t\le t'} S_\mathcal{C}(t)>s_0.
\]

Let us estimate the difference $\Delta S=S_\mathcal{C}-S_\Gamma^0$ on the 
interval 
$\varepsilon\le I\le i_0$.
From 
\[
\frac{dS_\mathcal{C}}{dI}=\frac{- R_0^\mathcal{C}(I)S_\mathcal{C}I+(1-S_\mathcal{C})\rho}{R^\mathcal{C}_0(I) S_\mathcal{C}I-I}, \qquad\frac{dS^0_\Gamma}{dI}=\frac{-R^\Gamma_0(I)S^0_\Gamma I+(1-S^0_\Gamma)\rho}{R^\Gamma_0(I) S^0_\Gamma I-I},
\]
after simple manipulations one obtains
\[
\frac{d}{dI} \Delta S= a(I)\Delta S+b(I),
\]
where
\[
a(I)=\frac{-\rho}{I(R^\mathcal{C}_0(I) S^0_\mathcal{C}(I)-1)(R^\Gamma_0(I) S^0_\Gamma(I)-1)}
\!\!\left(R^\Gamma_0(I) S^0_\Gamma(I)-1+R^\mathcal{C}_0(I)(1-S^0_\Gamma(I))-\frac{R^\mathcal{C}_0(I) I}{\rho}\right)\!\!,
\]
\[
b(I)=S^0_\Gamma(I)\bigl(R^\mathcal{C}_0(I)-R^\Gamma_0(I)\bigr)\left(1-\frac{ (1-S^0_\Gamma(I))\rho}{I}\right).
\]
Combining \eqref{hh} and the similar estimates $R^\Gamma_0(I) S^0_\Gamma(I)\ge R_0^{int}s_0>1$ with $R^\mathcal{C}_0(I)\le R_0^{nat}$ (cf.~\eqref{rrr}), and taking into account \eqref{bb3}, we see that $a(I)<0$ for $0<I\le i_0$, hence
\[
\frac{d}{dI} |\Delta S|= a(I)|\Delta S|+ b(I)\, \text{sign}\, \Delta S  \le |b(I)|\le \frac{|R^\mathcal{C}_0(I)-R^\Gamma_0(I)|}{I}
\]
and therefore
\[
|\Delta S(i_0)|\le |\Delta S(\varepsilon)|+\int_\varepsilon^{i_0}\frac{|R^\mathcal{C}_0(I)-R^\Gamma_0(I)|}{I}\,dI.
\]
From \eqref{LipP*}, it follows that 
\[
|R^\mathcal{C}_0(I)-R^\Gamma_0(I)|\le K \|r^{0,1}-r^{0,2}\|_{L_1(\Pi;\mathbb{R})}, \qquad \varepsilon\le I\le i_0,
\]
where $r^{0,1}$ (respectively, $r^{0,2}$) is the state functions of the Preisach operator at the moment when $I=\varepsilon$
for the trajectory $\mathcal C$ (respectively, $\Gamma_0$). But the compatibility condition \eqref{compatibility} implies that
$\|r^{0,1}-r^{0,2}\|_{L_1(\Pi;\mathbb{R})}\le\varepsilon$, hence
\[
|\Delta S(i_0)|\le |\Delta S(\varepsilon)|+\int_\varepsilon^{i_0}\frac{K\varepsilon}{I}\,dI \le\varepsilon+K\varepsilon\ln\frac{i_0}{\varepsilon}.
\]
We see that the point of intersection of the trajectory $\mathcal C$ with the line $I=i_0$ (which is independent of $\varepsilon$) converges 
to the trajectory $\Gamma_0$ as $\varepsilon\to 0$, and the convergence is uniform with respect to the initial conditions from the set $U_\varepsilon$ and all admissible initial states of the Preisach operator. Simultaneously, the state functions
of the Preisach operator corresponding to $I=i_0$ for the trajectories $\mathcal C$ and $\Gamma_0$ satisfy
$\|r^{1}-r^{2}\|_{L_1(\Pi;\mathbb{R})}\le\varepsilon$.
Therefore, the continuous dependence of trajectories on the initial point in the phase space $\mathfrak U$ of system \eqref{r0model}, \eqref{pre} implies that the point of intersection of $\mathcal C$ with the line $S=1/\hat R_0(I)$
converges to $M_0$ as $\varepsilon\to0$ for $\mathcal C(0)\in U_\varepsilon$, which together with \eqref{III} proves the theorem.

\end{document}